\documentclass[lettersize,journal]{IEEEtran}
\usepackage{amsmath,amsfonts}
\usepackage{algorithmic}
\usepackage{algorithm}
\usepackage{array}
\usepackage[caption=false,font=normalsize,labelfont=sf,textfont=sf]{subfig}
\usepackage{textcomp}
\usepackage{stfloats}
\usepackage{url}
\usepackage{verbatim}
\usepackage{graphicx}
\usepackage{cite}

\usepackage{booktabs}
\usepackage{multirow}
\usepackage{multicol}

\usepackage{pifont}  
\usepackage{orcidlink}

\usepackage{amsthm}     
\usepackage{tcolorbox}
\tcbuselibrary{theorems}

\newtheorem{definition}{Definition}[subsection]

\tcolorboxenvironment{theorem}{
  boxrule=0.5pt,
  boxsep=0pt,
  left=5pt,
  right=5pt,
  top=5pt,
  bottom=5pt,
  colback=white 
}

\tcolorboxenvironment{definition}{
  boxrule=0.5pt,
  boxsep=0pt,
  left=5pt,
  right=5pt,
  top=5pt,
  bottom=5pt,
  colback=white 
}

\begin{document}

\title{DRP-FLR: Data-Driven Assessment of Demand Response Potential for Flexible Load Regulation in Smart Grids}

\author{Yunhao~Yao~\orcidlink{0000-0002-7433-3145}, Siyu~Jing, Yang~Yang, Qiang~Xu, Changqi~Weng and Xiang-Yang~Li~\orcidlink{0000-0002-6070-6625},~\IEEEmembership{Fellow,~IEEE}
\thanks{Yunhao Yao, Siyu Jing, and Xiang-Yang Li are with the University of Science and Technology of China (USTC), China.}
\thanks{Yunhao Yao, Yang Yang, Qiang Xu, and Changqi Weng are with Anhui Zhongxin Jiyuan Information Technology Co., Ltd., China.}
\thanks{E-mail: yaoyunhao@mail.ustc.edu.cn, siyujing@mail.ustc.edu.cn, yang yang@zx-jy.com, xuqiang@zx-jy.com, wengcq@zx-jy.com, xiangyangli@ ustc.edu.cn.}
\thanks{Xiang-Yang Li is the corresponding author.}
\thanks{Manuscript received April 19, 2021; revised August 16, 2021.}}

\markboth{Journal of \LaTeX\ Class Files,~Vol.~14, No.~8, August~2021}%
{Shell \MakeLowercase{\textit{et al.}}: A Sample Article Using IEEEtran.cls for IEEE Journals}


\maketitle

\begin{abstract}
The rapid growth of AI workloads and renewable energy resources exacerbates supply-demand imbalance in power systems, making traditional load regulation designed for efficient allocation inadequate and motivating demand response (DR) mechanisms to enable load controllability in smart grids.
However, existing DR-oriented approaches either focus on optimizing electricity cost or occupant comfort with limited benefit to system-level balance. 
Others overlook the diverse and dynamic consumption patterns of heterogeneous energy entities, leading to significant over- or under-regulation.
Therefore, we propose DRP-FLR. 
First, DRP-FLR achieves accurate short-term load forecasting by embedding exogenous knowledge (e.g., entity information, prediction time) into historical load representations.
Next, it constructs entity-specific load-pattern profiles by clustering historical load curves, and estimates DR potential by matching forecasted loads with pattern profiles. 
Finally, DRP-FLR formulates flexible load regulation as a mixed-integer optimization problem and solves it with an MILP solver to jointly optimize DR utilization, participant economic benefit, and renewable accommodation, while enforcing supply-demand balance and economic feasibility.
Experiments on a regional grid and a campus microgrid show that DRP-FLR reduces regulation deviation by \(36.63\%\)-\(91.87\%\) and improves participant benefit by \(44.66\%\) on average.
\end{abstract}


\section{Introduction}
\label{intro}

The rapid advancement of artificial intelligence (AI) and the large-scale integration of renewable energy resources are reshaping modern power systems~\cite{stanelyte2022overview, mahmood2024impacts, wang2024ai}. 
On the demand side, AI-related computing workloads introduce fast-growing and bursty electricity consumption. 
On the supply side, renewables such as wind and solar bring significant variability and uncertainty. 
Together, these trends aggravate supply-demand imbalance and increase the operational difficulty of maintaining grid stability. 
In this context, smart grids have emerged as a key paradigm to enhance controllability through data-driven monitoring and automated decision-making. 
Among the available levers, \emph{demand response} mechanism enables flexible load regulation to support reliable grid operation by adjusting consumption in response to grid conditions~\cite{stanelyte2022overview, mahmood2024impacts}.

\textbf{Limitations of Traditional Load Regulation.} 
Traditional load regulation typically assumes that generation capacity is sufficient to meet demand, and they focus on efficient allocation under a supply-demand balance constraint~\cite{bashyal2025multi, alipour2017minlp, roh2015residential}. With the introduction of high-energy-demand AI workloads and volatile renewable generation, the grid may frequently face conditions where available generation is lower than demand. In such cases, purely supply-oriented or balance-assuming strategies may fail to provide effective control actions.

\textbf{Limitations of Existing DR-oriented Regulation.} 
Recent DR-oriented techniques also remain insufficient in practice. 
One line of work is designed for energy consumers (e.g., households), aiming to reduce electricity bills or improve comfort, which does not address the system-level supply-demand imbalance~\cite{eshraghi2019enhanced, veras2018multi}. 
Another line of work targets grid operators and makes peak shaving/valley filling decisions~\cite{li2024online, ruiz2019integration, muqtadir2025day, ruiz2024applications}. 
Yet many of them rely on coarse assumptions about load flexibility and do not account for the diverse consumption patterns across heterogeneous energy entities, which leads to severe over- or under-regulation, ultimately degrading both operational reliability and participant incentives.

These gaps motivate a load regulation mechanism built upon \emph{high-accuracy DR potential assessment}. 
To enable such a mechanism, several key challenges must be addressed:

\noindent~\underline{Accurate Future Load Forecasting.}
Both the discovery of regulation needs and the design of regulation strategies depend on accurate load forecasts. 
Future load is jointly affected by historical trends, time-of-day effects, day-of-week/holiday, and other exogenous factors. 
How to effectively fuse such multi-source knowledge into a forecasting model to achieve high-precision predictions remains a critical challenge.

\noindent~\underline{Adaptive DR Potential Estimation.}
A smart grid involves a large number of heterogeneous entities, each exhibiting multiple and dynamic consumption patterns. 
Designing a unified framework, which adapts to diverse entity types and their evolving load behaviours, can accurately estimate an entity's adjustable load range at a specific future time, and is essential for dependable DR-oriented regulation.

\noindent~\underline{Efficient Regulation Strategy Solving.}
Flexible load regulation is inherently a multi-objective optimization problem that jointly considers grid stability, economic feasibility, participant incentives, and renewable energy accommodation. 
This formulation often becomes a large-scale integer optimization problem, rendering brute-force search computationally intractable due to exponential complexity. 
Therefore, efficiently finding a regulation plan that satisfies operating constraints is essential.

To address these challenges, we propose \textbf{DRP-FLR}. 
To our knowledge, DRP-FLR is the first load regulation strategy to assess DR potential through explicit consumption-pattern analysis while jointly optimizing supply–demand balance, participant incentives, economic feasibility, and renewable energy accommodation.
Our main contributions are as follows:

\begin{itemize}
\item We design a forecast model that embeds multiple exogenous knowledge, including energy entity information and prediction time, into historical load representations, enabling more accurate short-term load forecasting.

\item We propose a pattern-driven method to estimate DR potential (i.e., the feasible regulation range) through consumption-pattern clustering. 
By recognizing consumption patterns, our approach improves the accuracy of DR potential assessment for heterogeneous entities.

\item We formulate a regulation model that maximizes DR potential utilization, participant economic benefits, and renewable energy accommodation, while ensuring supply-demand balance and economic feasibility. 
We further implement an efficient and accurate solver based on mixed-integer linear programming (MILP).

\item We implement DRP-FLR prototypes in both a regional grid and a campus microgrid, across various regulation scales and data granularities. 
Compared with existing DR-oriented regulation strategies, DRP-FLR substantially reduces regulation deviation by 36.63\%-91.87\%, while increasing average participant benefit by 44.66\%.
\end{itemize}

\section{Related Works}
\label{related works}

\subsection{Short-Term Load Forecast}
Short-term load forecasting (STLF) refers to predicting electricity demand over a short time horizon (e.g., minutes to several days) using historical load data and exogenous variables such as calendar information and weather conditions.
In load regulation, operators must anticipate the future states of a smart grid, and STLF provides the forward-looking load information required for scheduling and control decisions.

Statistical methods characterize load as a stochastic sequence and forecast future demand by fitting structured temporal patterns (e.g., trend and seasonality)~\cite{moslemi2024comprehensive, kim2023time, neshat2018nonlinear}.
Traditional machine learning methods formulate STLF as a regression problem that maps engineered features (e.g., historical lags and calendar variables) to future loads~\cite{chen2019short, raju2022approach}.
Recurrent Neural Network-based methods (e.g., LSTM/GRU) perform sequence-to-sequence learning leveraging hidden states to model temporal dependencies in load trends and can be extended to probabilistic forecasting~\cite{ochoa2025enhancing, nygaard2025enhancing, rahman2018predicting, guo2021short}.
Transformer-based methods leverage self-attention to capture long-range temporal dependencies and cross-variable interactions, often delivering strong performance under complex load patterns~\cite{chan2024transformer, zhao2021short}.
Hybrid frameworks combine heterogeneous forecast modules to improve accuracy and robustness under distribution shifts and diverse operating statuses~\cite{balakrishnan2025stacked, chow2021short, xiao2018hybrid, jnr2021hybrid}.

\subsection{Traditional Load Regulation}
Traditional load regulation typically studies how a power system allocates loads under the available supply by shaping the consumption behaviours of heterogeneous energy users, so as to maximize overall objectives such as production efficiency, energy cost savings, and residential quality of life.

One representative work models energy consumption and industrial production dynamics as a Partially Observable Markov Decision Process and leverages reinforcement learning to balance production efficiency against energy cost~\cite{bashyal2025multi}.
Another considers an energy hub with multiple energy inputs/outputs and proposes a $2m+1$ point-estimation-based probabilistic scheduling method to achieve supply-demand balance~\cite{alipour2017minlp}.
A residential scheduling study categorizes household appliances based on their usage patterns and optimizes each appliance's operating time and energy level to maximize user satisfaction while satisfying a budget constraint~\cite{roh2015residential}.

However, with increasing penetration of distributed renewables and the surge of urban energy demand, the regulation problem is more dominated by uncertainty and persistent scarcity rather than marginal cost optimization. 
Consequently, traditional load regulation is ineffective when generation capacity falls below demand.

\subsection{DR-Oriented Load Regulation}
With the increasing integration of distributed renewables and mismatch between generation capacity and demand, Demand Response has become a key mechanism to mitigate supply-demand imbalance through load-side flexibility~\cite{stanelyte2022overview, mahmood2024impacts}.
In a DR program, the grid operator publishes regulation requirements and operational constraints. 
Participants submit their regulation capacities, which the operator evaluates through historical performance to identify eligible customers for load adjustments that restore balance.

\subsubsection{Strategies Minimizing Consumer Costs}
Early DR-oriented load regulation strategies primarily aimed to minimize energy consumer costs.
In integrated cooling–heating–power systems with coupled storage and wind/PV integration, Eshraghi et al. coordinate the control of multiple storage devices to reduce electricity purchases~\cite{eshraghi2019enhanced}.
Under multi-horizon planning and real-time pricing, Veras et al. optimize appliance-level scheduling across device categories to balance user satisfaction and electricity cost~\cite{veras2018multi}.
However, consumer-oriented strategies cannot resolve the fundamental grid-level challenge of supply-demand mismatch. 
This has driven the development of operator-oriented strategies that select participants to deliver the required system-wide load reduction.

\subsubsection{Strategies Balancing Supply and Demand}
When generation capacity is insufficient, balancing supply and demand requires encouraging consumers to participate in DR for peak shaving, which in turn necessitates assessing each participant's DR potential for optimal selection.
Li et al. forecast load consumption with and without DR participation, thereby inferring the feasible regulation range~\cite{li2024online}.
Ruiz et al. derive adjustable load ranges from government statistics as a coarse proxy~\cite{ruiz2019integration}.
Muqtadir et al. approximate regulatable load mainly as air-conditioning demand and model it as positively correlated with ambient temperature~\cite{muqtadir2025day}.
In a separate study, Ruiz et al. apply quantile regression forests to predict load variability intervals, treating the interval width as the load adjustability~\cite{ruiz2024applications}.

Despite their effectiveness in some settings, these DR potential assessment methods have notable limitations: 
Estimating loads with and without DR participation requires fine-grained labels impractical to obtain, while government statistics can mask heterogeneity across consumers. 
Although some methods distinguish among different consumers, they still fail to capture diverse and evolving consumption patterns, which can lead to biased estimation.

\section{Problem Setup}
\label{prob}

\begin{figure}[t]
  \centering
  \includegraphics[width=1.0\columnwidth]{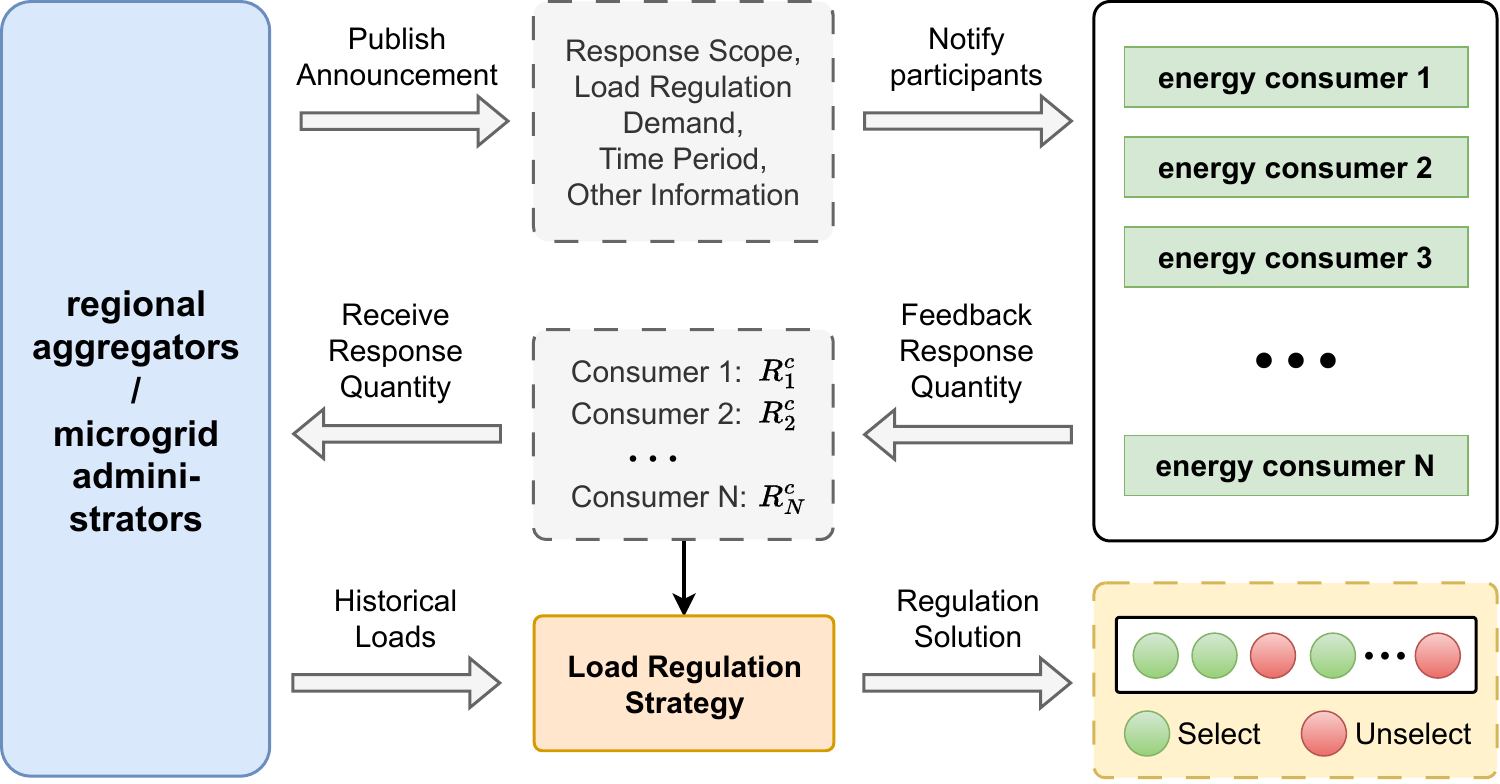}
  \caption{The Process of Demand Response Mechanism}
  \label{fig:DR}
\end{figure}

\subsection{Demand Response Mechanism}
Demand response (DR) is an essential mechanism within smart grids that enables dynamic adjustment of electricity consumption in response to time-varying grid conditions~\cite{stanelyte2022overview, mahmood2024impacts}. 
By guiding consumers to shift their power usage during peak periods to off-peak hours, DR helps balance power supply and demand to enhance grid stability.
In practice, grid operators (e.g., regional aggregators or microgrid administrators) identify potential peak loads based on load forecasts and assess whether available energy capacity can accommodate the predicted demand. 
When a supply-demand imbalance is anticipated, the DR mechanism is activated.

As illustrated in Fig.~\ref{fig:DR}, the DR process involves interactions between the grid operator and participating energy consumers. 
Initially, the operator issues a public announcement specifying the regulation demand, scheduling time and other relevant details. 
This information is delivered to potential participants, who evaluate their load flexibility and submit contractual response quantities $R_k^{C}$, representing their regulation capacities. 
The operator then aggregates these responses and executes a \textit{load regulation strategy} to determine the optimal regulation solution, including which consumers are selected for participation (\textit{Select}) and which are excluded (\textit{Unselect}). 
Finally, the selected consumers perform the corresponding load adjustments, thereby achieving real-time balance between supply and demand, improving grid reliability, and enhancing the overall utilization of renewable energy resources.

\subsection{Optimization Objective}
Our objective is to implement optimal load regulation for the smart grid within a DR framework. 
We consider a representative regional grid (e.g., a district within a city) or a microgrid (e.g., an industrial campus), 
comprising $N$ energy-consumption entities.
As a grid operator, our goal is to determine an optimal solution that simultaneously ensures grid stability, enhances economic efficiency, and promotes renewable energy integration. 
Accordingly, the optimization problem is formulated as follows:

\[
\max_{\mathbf{S}} \quad \sum_{k=1}^{N} S_k \cdot P_k \cdot W_k 
+ \alpha \cdot \sum_{k=1}^{N} S_k \cdot B_k \cdot W_k,
\]

\[
\begin{aligned}
\text{s.t.} \quad 
& L \leq \sum_{k=1}^{N} S_k \cdot R_k \leq (1 + \sigma) \cdot L, \\
& S_k \in \{0, 1\}, \quad \forall k \in \{1, 2, \ldots, N\}.
\end{aligned}
\]

Here, 
$\mathbf{S} = [S_1, S_2, \ldots, S_N] \in \{0,1\}^{1 \times N}$ denotes the selection vector, where $S_k = 1$ indicates that the $k$-th entity is selected to participate in load regulation.
$\mathbf{P} = [P_1, P_2, \ldots, P_N]$ represents the DR potential of each entity, reflecting its ability to meet or exceed the contracted regulation volume. 
Selecting entities with greater potential enhances the probability of maintaining the supply-demand balance.  
$\mathbf{B} = [B_1, B_2, \ldots, B_N]$ denotes the economic benefits obtainable by each participant, where higher benefits incentivize more active participation.
$\mathbf{W} = [W_1, W_2, \ldots, W_N]$ represents weighting factors, where renewable energy participants are assigned higher weights to encourage their involvement.
$\mathbf{R} = [R_1, R_2, \ldots, R_N]$ indicates the load regulation capacity of each entity.
The total dispatched regulation must satisfy the regulation demand $L$ without exceeding $(1 + \sigma) \cdot L$, ensuring both stability and economic feasibility.

\begin{figure*}[t]
  \centering
  \includegraphics[width=1.6\columnwidth]{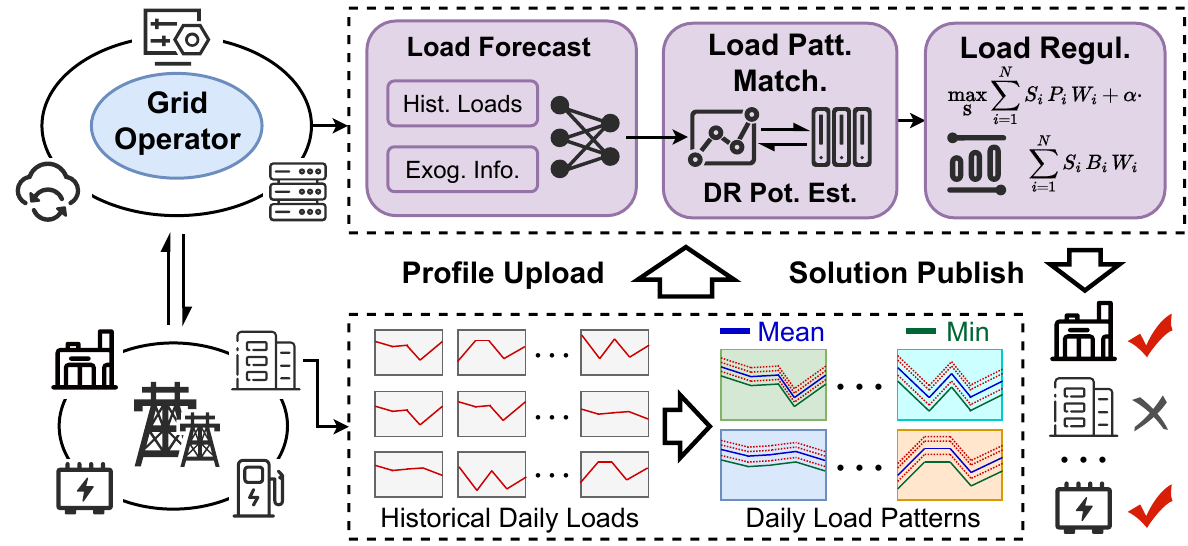}
  \caption{The System Overview of DRP-FLR}
  \label{fig:system}
\end{figure*}

\section{System Design}
\label{sys}

In this section, we provide a detailed explanation of DRP-FLR. 
First, the \texttt{Load Forecast} module predicts each energy consumer's future demand to support load-pattern matching. 
Based on the matched consumption patterns, the \texttt{DR Potential Estimation} module quantifies each consumer's potential range of load reduction.
Finally, the \texttt{Load Regulation} module selects an optimal solution that balances grid stability, economic efficiency, and renewable energy integration.
The system overview is shown in Fig.~\ref{fig:system}

\subsection{Multi-Knowledge-Embedded Future Load Forecast}

Since load forecasting is inherently a time-series prediction task, we introduce TimeGPT~\cite{garza2023timegpt} as the backbone predictor, denoted as $F_{\mathrm{pre}}(\cdot)$, which is a Transformer-based model.
It typically consists of:
(i) an embedding module that projects raw observations into latent tokens;
(ii) stacked Transformer blocks with multi-head self-attention and feed-forward layers to capture long-range temporal dependencies;
and (iii) a forecast head that produces one- or multi-step predictions.

Electricity consumption is not solely determined by historical loads; it is also influenced by exogenous knowledge related to the energy consumer and the calendar context.
For example, different types of enterprises exhibit different consumption patterns (e.g., manufacturing plants versus office buildings),
and the load profiles on workdays can differ substantially from those on holidays (e.g., shifted peak hours and increased operating loads).
Therefore, incorporating structured exogenous knowledge is crucial for robust and accurate forecasting.

Therefore, we propose a multi-knowledge-embedded forecasting method.
Specifically, we encode the \emph{energy-consumer type}, \emph{forecasting time period}, and \emph{workday/holiday indicator} into learnable embeddings,
and overlay them onto the embeddings of the historical load sequence.
The fused embeddings are fed into the backbone $F_{\mathrm{pre}}(\cdot)$ for hourly load prediction.

Let $\mathrm{Ent}^k$, $\mathrm{Time}$, and $\mathrm{DayType}$ denote the exogenous knowledge, and let
$\mathrm{Emb}^k_{\mathrm{ent}}$, $\mathrm{Emb}_{\mathrm{time}}$ and $\mathrm{Emb}_{\mathrm{dayType}}$ be their embeddings.
Let $\mathrm{Load}^k[M-\tau:M]$ denote the historical load window and $\mathrm{Emb}^k_{\mathrm{load}}[M-\tau:M]$ its embedding sequence.
The full forecasting pipeline is formulated as:
\[
\begin{aligned}
& \mathrm{Emb}^k_{\mathrm{load}}[M-\tau:M] = \phi_{\mathrm{load}}\!\left(\mathrm{Load}^k[M-\tau:M]\right), \\
& \mathrm{Emb}^k_{\mathrm{ent}} = \phi_{\mathrm{ent}}(\mathrm{Ent}^k), \\
& \mathrm{Emb}_{\mathrm{time}} = \phi_{\mathrm{time}}(\mathrm{Time}), \\
& \mathrm{Emb}_{\mathrm{dayType}} = \phi_{\mathrm{dayType}}(\mathrm{DayType}), \\
& \mathrm{Emb}^k_{\mathrm{exo}} = \mathrm{Emb}^k_{\mathrm{ent}} + \mathrm{Emb}_{\mathrm{time}} + \mathrm{Emb}_{\mathrm{dayType}}, \\
& \widetilde{\mathrm{Emb}}^k[M-\tau:M] = \mathrm{Emb}^k_{\mathrm{load}}[M-\tau:M] \oplus \mathrm{Emb}^k_{\mathrm{exo}}, \\
& \mathrm{Load}^k_{\mathrm{pre}}[M+1] = F_{\mathrm{pre}}\!\left(\widetilde{\mathrm{Emb}}^k[M-\tau:M]\right).
\end{aligned}
\]

\subsection{Consumption-Pattern-Driven DR Potential Estimation}

\subsubsection{Energy Consumption Pattern Analysis}

Different energy consumers (e.g., shopping malls, factories, and residential communities) exhibit distinct energy consumption patterns due to their heterogeneous operating schedules, occupancy profiles, and production activities.
Moreover, even for the same energy consumer, its consumption pattern can vary over time under different operating conditions (e.g., seasonal changes, special events, or urgent production demands).

Different consumption patterns lead to different \emph{load regulation} capabilities.
For instance, a residential community may have a limited regulation range under routine daily life conditions.
However, this range can change substantially under atypical conditions (e.g., abrupt weather shifts) or when the community deploys additional flexible resources.
Similarly, industrial consumers may present substantially different regulation ranges under normal versus order-rush production states, where the priority shifts to maintaining throughput, thereby narrowing the feasible controllable load range.

To estimate the load regulation capacity $R_k$ $(1 \leq k \leq N)$ of each energy-consumption entity, a prerequisite is to identify and analyze the \emph{energy consumption patterns} that the entity operates in, because the controllable ranges are pattern-dependent.
In practice, the number of patterns for each energy consumer is unknown \emph{in advance}.
Therefore, we adopt the density-based clustering algorithm DBSCAN~\cite{ester1996density} to cluster the historical \emph{daily load curves} of $N$ energy consumers and thereby identify their consumption patterns.

\begin{figure*}[t]
  \centering
  \includegraphics[width=1.8\columnwidth]{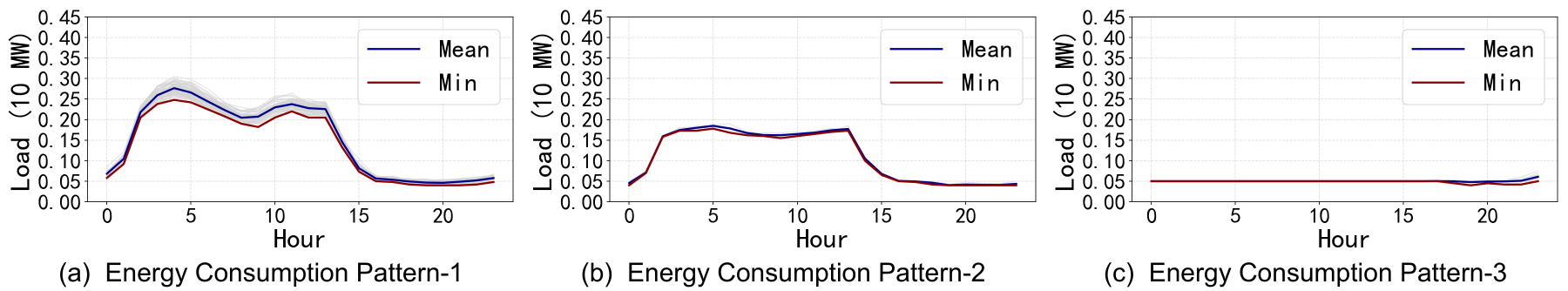}
  \caption{Example of Demand Response Potential: An Estate Company in A Regional Grid}
  \label{fig:ex_pattern_regional_grid}
\end{figure*}

\begin{figure*}[t]
  \centering
  \includegraphics[width=1.8\columnwidth]{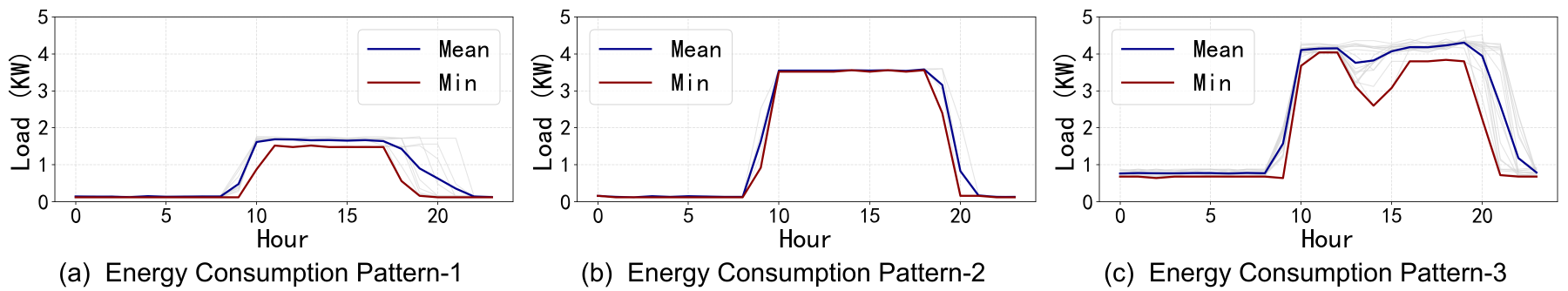}
  \caption{Example of Demand Response Potential: An Office Area in A Microgrid}
  \label{fig:ex_pattern_micro_grid}
\end{figure*}

(1) Each daily load curve contains 24 points, where each point corresponds to the hourly average load of that day.
Let $\mathbf{L}^k$ denote the collection of historical daily load curves for the $k$-th energy entity, and let $\mathbf{L}^k_i \in \mathbb{R}^{24}$ denote the daily load curve of entity $k$ on day $i$.
Suppose that entity $k$ has $D_k$ historical days, then $\mathbf{L}^k = \{\mathbf{L}^k_i\}_{i=1}^{D_k}$.

(2) DBSCAN is chosen because it does not require pre-specifying the number of clusters and can explicitly identify outliers (i.e., noise days) caused by abnormal operations or data-quality issues.
To select appropriate DBSCAN hyperparameters, we first compute the pairwise Euclidean distances between daily load curves.
For any two days $i$ and $j$ of entity $k$, the Euclidean distance is $dis^k_{i,j} = \lVert \mathbf{L}^k_i - \mathbf{L}^k_j \rVert_2$.

(3) For each daily load sample, we further compute the mean distance to its four nearest neighbours within the same entity's historical load collection.
Let $\mathcal{N}^k_4(i)$ denote the index set of the four nearest neighbors of $\mathbf{L}^k_i$ within $\{\mathbf{L}^k_j\}_{j\neq i}$ with respect to $dis^k_{i,j}$.
We define the mean 4-NN distance as: 
\[
\bar{dis}^k_i = \frac{1}{4}\sum_{j \in \mathcal{N}^k_4(i)} dis^k_{i,j}.
\]

(4) We then take the 25\%, 50\%, and 75\% quantiles of these mean nearest-neighbour distances as candidate values for the neighbourhood radius $\varepsilon$.
We collect the mean 4-NN distance of each daily load curve for entity $k$ as:
\[
\bar{\mathcal{D}}^k = \{\bar{dis}^k_i\}_{i=1}^{D_k}.
\]
Let $Q_p(\cdot)$ denote the $p$-quantile operator. The candidate set of neighbourhood radius is:
\[
\mathcal{E}^k =
\left\{
Q_{0.25}\!\left(\bar{\mathcal{D}}^k\right),
Q_{0.50}\!\left(\bar{\mathcal{D}}^k\right),
Q_{0.75}\!\left(\bar{\mathcal{D}}^k\right)
\right\}.
\]

(5) Meanwhile, we set $\mathrm{minPts} \in \{3, 4\}$.
For each candidate pair $(\varepsilon,\mathrm{minPts})$, we run DBSCAN and select the configuration that maximizes the silhouette coefficient~\cite{rousseeuw1987silhouettes} computed on the non-noise points.
Let $\mathcal{M} = \{3,4\}$ denote the candidate set for $\mathrm{minPts}$.
For each $(\varepsilon,m)\in \mathcal{E}^k \times \mathcal{M}$, we run DBSCAN and obtain cluster labels:
\[
\mathbf{c}^k(\varepsilon,m)=\mathrm{DBSCAN}\!\left(\mathbf{L}^k;\varepsilon,m\right),
\]
where noise points are assigned the label $-1$.

(6) Finally, we choose the hyperparameter pair that achieves the highest silhouette score, yielding the optimal clustering configuration.
Let $\mathcal{I}^k_{\mathrm{core}}(\varepsilon,m)=\{i \mid \mathbf{c}^k_i(\varepsilon,m)\neq -1\}$ denote the set of non-noise days.
We compute the silhouette coefficient using only the non-noise points:
\[
\begin{aligned}
\mathrm{Score}^k(\varepsilon,m)
&= \mathrm{Silhouette}\!\left(
\{\mathbf{L}^k_i\}_{i\in \mathcal{I}^k_{\mathrm{core}}(\varepsilon,m)},
\right.\\
&\qquad \qquad \qquad \left.
\{\mathbf{c}^k_i(\varepsilon,m)\}_{i\in \mathcal{I}^k_{\mathrm{core}}(\varepsilon,m)}
\right).
\end{aligned}
\]
The optimal hyperparameters for entity $k$ are then given by:
\[
(\varepsilon_k^{\ast}, m_k^{\ast})
=
\arg\max_{(\varepsilon,m)\in \mathcal{E}^k \times \mathcal{M}}
Score^k(\varepsilon,m),
\]
and the corresponding optimal patterns are
$
\mathbf{c}^k(\varepsilon_k^{\ast}, m_k^{\ast}).
$

\subsubsection{Demand Response Potential Assessment}
For the $k$-th energy-consumption entity, its hourly demand response potential is determined by two variables:
(i) the contracted response amount $R^c_k$ pre-signed between the entity and the grid operator,
and (ii) the entity's load regulation capacity $R_k$ at the target time.
Intuitively, $R^c_k$ reflects the \emph{committed} DR capability, while $R_k$ reflects the \emph{feasible} regulation range inferred from the predicted operating pattern.

Assume the current time is $M$, and the intended load regulation time is $M+\delta_T$, where $1 \leq \delta_T \leq 24$.
We first predict the next 24-hour load curve for the $k$-th entity as:
$\mathrm{Load}^k_{\mathrm{pre}}[M+1] = F_{\mathrm{pre}}\!\left(\widetilde{\mathrm{Emb}}^k[M-\tau:M]\right)$,
where $\widetilde{\mathrm{Emb}}^k[M-\tau:M]$ denotes the fused (load + exogenous knowledge) embedding sequence over the look-back window.
More generally, we denote the predicted day-ahead profile as:
\[
\mathbf{PL}^k_M = [
\mathrm{Load}^k_{\mathrm{pre}}[M+1],\,
\dots,\, \mathrm{Load}^k_{\mathrm{pre}}[M+24]
]^{\top} \in \mathbb{R}^{24}.
\]

According to \textbf{Section IV.B.1}, $\mathbf{c}^k(\varepsilon_k^{\ast}, m_k^{\ast})$ is the optimal daily-load-pattern clustering for the $k$-th entity obtained by DBSCAN,
where $\varepsilon_k^{\ast}$ and $m_k^{\ast}$ are the optimal neighborhood radius and $\mathrm{minPts}$, respectively.
We assign the predicted day-ahead profile $\mathbf{PL}^k_M$ to its corresponding pattern cluster:
\[
Grp^k_M =
\mathrm{Assign}\!\left(
\mathbf{PL}^k_M,\,
\mathbf{c}^k(\varepsilon_k^{\ast}, m_k^{\ast})
\right),
\]
where $Grp^k_M$ is the cluster (pattern) index. If $Grp^k_M$ is identified as noise, we map it to the nearest non-noise cluster.

Let $\mathcal{C}(Grp^k_M)$ denote the set of historical days of the $k$-th entity that belong to cluster $Grp^k_M$.
For any hour offset $\delta_T$, we define the regulation capacity $R_k(M+\delta_T)$ as the difference between the \emph{mean} and the \emph{minimum} load within the identified pattern cluster:
\[
\mu_{Grp^k_M}(\delta_T)
=
\frac{1}{|\mathcal{C}(Grp^k_M)|}
\sum_{i \in \mathcal{C}(Grp^k_M)}
\mathbf{L}^k_i[\delta_T],
\]
\[
\ell_{Grp^k_M}(\delta_T)
=
\min_{i \in \mathcal{C}(Grp^k_M)}
\mathbf{L}^k_i[\delta_T],
\]
\[
R_k(M+\delta_T)
=
\mu_{Grp^k_M}(\delta_T) - \ell_{Grp^k_M}(\delta_T),
\]
where $\mathbf{L}^k_i[\delta_T]$ is the load of entity $k$ on day $i$ at hour $\delta_T$.

Finally, we define the demand response potential of the $k$-th entity at time $M+\delta_T$ as the ratio between its estimated regulation capacity and its contracted response amount, as shown in \textbf{Definition~\ref{def:potential}}.
Figs.~\ref{fig:ex_pattern_regional_grid} and~\ref{fig:ex_pattern_micro_grid} illustrate representative examples of clustering energy-consumption patterns for entities in the regional grid and the microgrid, and of assessing their demand-response potential.

\begin{definition}[Demand Response Potential]
\label{def:potential}
For an energy-consumption entity $k$, given the contracted response amount $R^c_k(M+\delta_T)>0$ and the estimated load regulation capacity $R_k(M+\delta_T)$ at time $M+\delta_T$ $(1 \leq \delta_T \leq 24)$, the \emph{demand response potential} is defined as:
\[
P_k(M+\delta_T) = \frac{R_k(M+\delta_T)}{R^c_k}.
\]
\end{definition}

\subsection{Multi-Objective Load Regulation Solution}

\subsubsection{Optimization Objective Refinement}
According to the optimization objective in \textbf{Section III.B}, for the regulation time $M+\delta_T$, we aim to jointly maximize
(i) the aggregated DR potential of selected entities and
(ii) the aggregated expected benefit of participating entities, while encouraging renewable-friendly entities via weights.
Specifically, let $S_k(M+\delta_T)\in\{0,1\}$ be the binary selection variable indicating whether the $k$-th entity is scheduled to participate at $M+\delta_T$.
Let $P_k(M+\delta_T)$ denote the DR potential and $B_k(M+\delta_T)$ the predicted benefit.
The two objectives can be refined as:
\[
\max \; \sum_{k=1}^{N} S_k(M+\delta_T)\cdot P_k(M+\delta_T)\cdot W_k,
\]
\[
\max \; \sum_{k=1}^{N} S_k(M+\delta_T)\cdot B_k(M+\delta_T)\cdot W_k.
\]

By the definition of $P_k(M+\delta_T)$ in \textbf{Definition~\ref{def:potential}}, a larger potential indicates an entity with higher credibility (i.e., more reliable regulation capacity relative to its contracted amount) in demand response events.
Meanwhile, maximizing the predicted benefit encourages broader participation.
Following the demand response settlement rule of State Grid Corporation of China (SGCC)~\cite{NDRC2023DSM}, the expected benefit of the $k$-th entity is determined by its estimated regulation capacity $R_k(M+\delta_T)$ relative to its contracted response amount $R^c_k(M+\delta_T)$, as shown in \textbf{Definition~\ref{def:benefit}}.
Besides, the weight $W_k$ provides an explicit lever to promote renewable energy accommodation by favouring renewable-related entities in the regulation decision (e.g., assigning a higher $W_k$).

\begin{definition}[Demand Response Benefit]
\label{def:benefit}
Let $A$ denote the unit reward per regulated load amount, then $B_k(M+\delta_T)$ is defined piecewise as:
\[
t = M+\delta_T,
\]
\[
B_k(t)=
\begin{cases}
0, & R_k(t) < 0.8\,R_k^c(t), \\[4pt]
A \cdot R_k(t), & 0.8\,R_k^c(t)\le R_k(t)\le 1.2\,R_k^c(t), \\[4pt]
A \cdot 1.2 \cdot R_k^c(t), & R_k(t) > 1.2\,R_k^c(t).
\end{cases}
\]
\end{definition}

According to the constraints in \textbf{Section III.B}, the total regulation capacity at time $M+\delta_T$ should satisfy a bounded requirement:
\[
\begin{aligned}
    L(M+\delta_T) & \leq \sum_{k=1}^{N} S_k(M+\delta_T)\cdot R_k(M+\delta_T) \\
    & \leq (1+\sigma)\cdot L(M+\delta_T),
\end{aligned}
\]
where $L(M+\delta_T)$ is the required regulation demand at time $M+\delta_T$ and $\sigma \geq 0$ controls the tolerated over-regulation margin.
This constraint is designed to balance (i) supply-demand adequacy (no less than $L(M+\delta_T)$) and (ii) economic feasibility (avoiding excessive regulation beyond $(1+\sigma) \cdot L(M+\delta_T)$).

\subsubsection{MILP Problem Solving}
With binary decision variables $S_k(M+\delta_T)\in\{0,1\}$ and linear objective/constraint forms, the above problem can be modelled as a 0-1 integer programming problem.
Therefore, we can solve it using a standard mixed-integer linear programming (MILP) solver.

We implement the optimization using PuLP and its default open-source backend solver \texttt{PULP\_CBC\_CMD}.
\texttt{PULP\_CBC \_CMD} is PuLP's built-in command-line interface to the CBC (COIN-OR Branch-and-Cut) solver, an open-source MILP solver from the COIN-OR project.
Given a MILP model, CBC typically proceeds as follows: it first solves the LP relaxation, then applies a \emph{branch-and-bound} search to enforce integrality of integer variables, and accelerates convergence by generating \emph{cutting planes} to tighten the relaxation.

In our load regulation task, the solver is provided with:
(i) decision variables (here $S_k(M+\delta_T)$ and auxiliary slack variables if needed),
(ii) a linear objective function (single-objective form after scalarization), and
(iii) a set of linear constraints.
And it returns
(i) feasibility status,
(ii) the optimal/near-optimal values of decision variables,
and (iii) the achieved objective value(s).
From these outputs, we obtain the selected entity set and the total regulation capacity.

\subsubsection{Feasibility-Robust Design}
In real deployments, the strict regulation bound
$
L \le \sum_{k=1}^{N} S_k(M+\delta_T)\,R_k(M+\delta_T) \le (1+\sigma) \cdot L
$
may become infeasible due to limited flexibility or prediction errors.
To guarantee solvability while still discouraging constraint violations, we introduce two penalty parameters:
the \emph{shortfall} $\xi(M+\delta_T) \geq 0$ for the lower-bound violation and the \emph{slack} variable $\zeta(M+\delta_T) \geq 0$ for the upper-bound violation.
Let
$
\mathrm{sum\_R}(M+\delta_T)=\sum_{k=1}^{N} S_k(M+\delta_T)\cdot R_k(M+\delta_T).
$
We replace the hard bounds with the following \emph{soft} constraints:
\[
\mathrm{sum\_R}(M+\delta_T) + \xi(M+\delta_T) \ge L,
\]
\[
\mathrm{sum\_R}(M+\delta_T) \le (1+\sigma)\,L + \zeta(M+\delta_T),
\]
where $\xi$ absorbs any unavoidable under-delivery below $L$, and $\zeta$ absorbs any unavoidable over-delivery beyond $(1+\sigma) \cdot L$.

To avoid frequent violations, we incorporate large penalties into the objective. The resulting single-objective MILP is:
\[
\begin{aligned}
\max \sum_{k=1}^{N} & S_k(M+\delta_T) \cdot
\Big(P_k(M+\delta_T) + \alpha\,B_k(M+\delta_T) \Big)\,W_k \\
& - \lambda_{\xi} \cdot \xi(M+\delta_T) - \lambda_{\zeta} \cdot \zeta(M+\delta_T),
\end{aligned}
\]
where $\lambda_{\xi} \gg 1$ and $\lambda_{\zeta}\gg 1$ are large penalty coefficients.
Typically, we set $\lambda_{\xi} \ge \lambda_{\zeta}$ to prioritize satisfying the minimum required regulation $L$ over avoiding mild over-regulation.

\section{Experiments}
\label{exp}

\subsection{Experimental Setup}

\subsubsection{Datasets}
Our experiments are conducted on two representative flexible load regulation scenarios in smart grids: (i) a regional grid with heterogeneous industrial/commercial customers, and (ii) a campus-level microgrid with fine-grained submeters and integrated renewable generation and storage.

\textbf{The Regional Grid Dataset} contains 15-minute load measurements from $441$ enterprises located in an administrative district of a Chinese city, spanning from October~2022 to October~2023.
The covered sectors include real estate, industry, commerce, manufacturing, and public services, reflecting a diverse energy consumption portfolio.
Among them, $20$ have installed renewable energy resources (e.g., photovoltaic and wind generation), enabling us to evaluate renewable-aware scheduling and weighting strategies in demand response.

\textbf{The Microgrid Dataset} is collected from a smart industrial campus managed by a Chinese enterprise.
It consists of 15-minute measurements from $73$ submeters over the period January~2025 to March~2026.
The metering scope covers typical end-use categories and distributed energy resources, including: lighting and HVAC loads at the floor level, elevator loads, server/IT room loads, electric vehicle charging stations, photovoltaic generation, wind generation, and the charge/discharge load of an energy storage system.

For both datasets, we resample the original 15-minute time series to an hourly resolution to align the data with the hour-ahead regulation horizon used in our experiments.
Missing values are estimated via mean-based interpolation.
We adopt a unified sign convention: electricity consumption is represented as positive values, while electricity generation is represented as negative values.
This convention allows net power trajectories to be modelled consistently when jointly considering demand, renewable generation, and storage operation.

\subsubsection{Evaluation Metrics}
We evaluate DRP-FLR from two complementary perspectives: (i) \emph{regulation fulfilment} and (ii) \emph{participant incentive}.
Accordingly, our primary evaluation metrics are \textbf{Achievement Rate} and \textbf{Participant Benefit}.

\textbf{Achievement Rate} measures how well the selected entities fulfill the required regulation demand.
Let $L$ denote the regulation demand and $S_k\in\{0,1\}$ indicate the selection of entity $k$.
Denote the \emph{real load reduction} delivered by entity $k$ as $R_k^{\mathrm{real}} \ge 0$.
Then the Achievement Rate is defined as:
\[
\mathrm{AR} \;=\;
\frac{\sum_{k=1}^{N} S_k \cdot R_k^{\mathrm{real}}}{L}.
\]
Ideally, $\mathrm{AR}$ should lie within the target interval $[1,\,1+\sigma]$, where $\sigma \ge 0$ is the tolerated over-regulation margin.
An $\mathrm{AR}<1$ indicates insufficient regulation, while $\mathrm{AR}>1+\sigma$ indicates excessive regulation beyond the economic tolerance.

$\mathrm{AR}$ reflects two capabilities of a load regulation solution:
\begin{enumerate}
  \item \textit{Accuracy of demand response potential estimation.}
  In our optimization, the term
  $\sum_{k=1}^{N} S_k \cdot P_k \cdot W_k$
  drives DRP-FLR to prioritize entities with high estimated DR potential $P_k$.
  If $P_k$ is inaccurately estimated, DRP-FLR may select suboptimal entities whose real reductions are lower or higher than expected, leading to a too small or too large Achievement Rate.

  \item \textit{Ability to satisfy grid stability and economic feasibility.}
  Our constraints explicitly require the selected portfolio to ensure supply-demand balance while avoiding excessive dispatch:
  $L \le \sum_{k=1}^{N} S_k \cdot R_k,
  \,
  \sum_{k=1}^{N} S_k \cdot R_k \le (1+\sigma)\cdot L.$
  When $\mathrm{AR}\in[1,\,1+\sigma]$, the real load reduction aligns well with these constraints, indicating that the solution
  provides adequate regulation for grid stability without incurring unnecessary costs.
\end{enumerate}

\textbf{Participant Benefit} quantifies the total monetary reward that selected entities can obtain under their real regulation performance.
Let $B_k^{\mathrm{real}}$ denote the real load regulation of entity $k$, computed according to \textbf{Definition~\ref{def:benefit}}.
Then Participant Benefit is defined as:
\[
\mathrm{PB} \;=\; \sum_{k=1}^{N} S_k \cdot B_k^{\mathrm{real}}.
\]

This metric is also closely related to the accuracy of DR potential estimation.
In the objective, the term 
$\sum_{k=1}^{N} S_k \cdot B_k \cdot W_k$
is computed from the estimated regulation capacity $R_k$,
which implicitly assumes that the real reduction will be close to $R_k$.
If the potential estimation is inaccurate, the gap between $R_k$ and $R_k^{\mathrm{real}}$ becomes large, and selected entities may fail to achieve the expected reward.
As a result, the real Participant Benefit decreases, indicating weaker incentives and reduced willingness to participate in future regulation events.

\begin{table*}[t]
    \centering
    \caption{Evaluation Achievement Rate and Participant Benefit (/10 MW) Across Regulation Ratios and Trade-off Factors: Regional Grid}
    \resizebox{1.0\linewidth}{!}{
        \begin{tabular}{c|c|c|c|c|c|c}
        \toprule
        \midrule
        Regulation Ratios & \multicolumn{2}{c|}{6.0\%} & \multicolumn{2}{c|}{8.0\%} & \multicolumn{2}{c}{10.0\%} \\
        \midrule
        Trade-off Factor & Achievement Rate & Participant Benefit & Achievement Rate & Participant Benefit & Achievement Rate & Participant Benefit  \\
        \midrule
        1e1 & 118.67\% & 0.0227 & 104.65\% & 0.0341 & 95.05\% & 0.0604 \\
        \midrule
        1e2 & 118.73\% & 0.0228 & 104.66\% & 0.0343 & 95.08\% & 0.0605 \\
        \midrule
        1e3 & 119.62\% & 0.0237 & 105.15\% & 0.0355 & 95.22\% & 0.0614 \\
        \midrule
        1e4 & 133.53\% & 0.0311 & 113.11\% & 0.0417 & 96.50\% & 0.0635 \\
        \midrule
        \bottomrule
        \end{tabular}
    }
    \label{tab:ablation_regional}
\end{table*}

\begin{table*}[t]
    \centering
    \caption{Evaluation Achievement Rate and Participant Benefit (/KW) Across Regulation Ratios and Trade-off Factors: Microgrid}
    \resizebox{1.0\linewidth}{!}{
        \begin{tabular}{c|c|c|c|c|c|c}
        \toprule
        \midrule
        Regulation Ratios & \multicolumn{2}{c|}{6.0\%} & \multicolumn{2}{c|}{8.0\%} & \multicolumn{2}{c}{10.0\%} \\
        \midrule
        Trade-off Factor & Achievement Rate & Participant Benefit & Achievement Rate & Participant Benefit & Achievement Rate & Participant Benefit  \\
        \midrule
        1e1 & 103.42\% & 0.1976 & 97.87\% & 0.3056 & 94.51\% & 0.4613 \\
        \midrule
        1e2 & 103.97\% & 0.2279 & 99.41\% & 0.3915 & 95.54\% & 0.6091 \\
        \midrule
        1e3 & 105.40\% & 0.2735 & 100.56\% & 0.4617 & 96.99\% & 0.6910 \\
        \midrule
        1e4 & 122.71\% & 0.4097 & 120.61\% & 0.6677 & 115.74\% & 0.9406 \\
        \midrule
        \bottomrule
        \end{tabular}
    }
    \label{tab:ablation_micro}
\end{table*}

\subsection{Ablation Study}

\subsubsection{Various Regulation Ratios and Trade-off Factors}
We conduct DR-oriented load regulation experiments in both the \emph{regional distribution grid} and the \emph{campus microgrid} scenarios.
To validate the effectiveness and robustness of DRP-FLR under different operating conditions, we vary
(i) the regulation demand level, implemented as different \emph{regulation ratios}, and
(ii) the regulation preference, implemented by the trade-off factor \(\alpha\) that balances \emph{response potential} versus \emph{participant benefit}.
Unless otherwise specified, we set the tolerance margin in the regulation constraints to \(\sigma=0.2\), meaning that we expect the realized dispatch outcome to fall within \(1\) to \(1.2\) times the target regulation demand.
The experimental results are reported in TABLE~\ref{tab:ablation_regional} and TABLE~\ref{tab:ablation_micro}.

In the regional distribution grid, the Achievement Rate ranges from \(95.05\%\) to \(133.53\%\),
corresponding to a deviation of \([-4.95\%,\,13.53\%]\) relative to the ideal interval \([100\%,\,120\%]\).
In the campus microgrid, the Achievement Rate ranges from \(94.51\%\) to \(122.71\%\),
corresponding to a deviation of \([-5.49\%,\,2.71\%]\) relative to \([100\%,\,120\%]\).
The main reason for these deviations is that our regulation decision relies on \emph{estimated} regulation capacities derived from \emph{forecasted} future loads.
Even if the predicted regulation lies in the desired interval, forecasting errors can lead to discrepancies between the predicted and real load reductions, causing the Achievement Rate to under- or over-shoot the target range.

As the regulation ratio increases, we observe a gradual decrease in the Achievement Rate.
This phenomenon is rooted in the discrete nature of demand response participation: for an entity \(k\), the real response is either not participating or delivering a certain reduction \(R_k^{\mathrm{real}}\), rather than being continuously adjustable to any level.
When the regulation ratio is small, it becomes harder for DRP-FLR to identify a subset of entities whose aggregated reduction
\(
\sum_{k=1}^{N} S_k \, R_k^{\mathrm{real}}
\)
matches demand \(L\) tightly.
Consequently, the solution tends to exhibit "spillover" (over-regulation), resulting in a higher Achievement Rate.
As \(L\) increases, the aggregation of more entities provides finer combinatorial granularity, and the spillover becomes smaller, leading to a reduced Achievement Rate.

As the trade-off factor \(\alpha\) increases, DRP-FLR places more emphasis on participant benefit and tends to select more entities. 
This increases both the total regulation and the aggregate reward, resulting in a higher Achievement Rate and Participant Benefit.
When \(\alpha\) becomes sufficiently large, DRP-FLR may even prefer solutions that exceed the upper bound \((1+\sigma) \cdot L\) and accept the corresponding penalty, in order to include additional participants and improve total benefit.

\begin{table}[t]
    \centering
    \caption{Achievement Rate and Participant Benefit (/10 MW) with and without Consumption Pattern Clustering: Regional Grid}
    \resizebox{1.0\linewidth}{!}{
    \begin{tabular}{c|c c|c c|c c}
        \toprule
        \midrule
        Ragulation Ratio & \multicolumn{2}{c|}{6.0\%} & \multicolumn{2}{c|}{8.0\%} & \multicolumn{2}{c}{10.0\%} \\
        \midrule
        Clustering? & \ding{51} & \ding{55} & \ding{51} & \ding{55} & \ding{51} & \ding{55} \\
        \midrule
        Achievement Rate & 119.62\% & 210.87\% & 105.15\% & 169.92\% & 95.22\% & 125.17\% \\
        \midrule
        Participant Benefit & 0.0237 & 0.0377 & 0.0355 & 0.0418 & 0.0614 & 0.0462 \\
        \midrule
        \bottomrule
    \end{tabular}
    }
    \label{tab:clustering_regional}
\end{table}

\begin{table}[t]
    \centering
    \caption{Achievement Rate and Participant Benefit (/KW) with and without Consumption Pattern Clustering: Microgrid}
    \resizebox{1.0\linewidth}{!}{
    \begin{tabular}{c|c c|c c|c c}
        \toprule
        \midrule
        Ragulation Ratio & \multicolumn{2}{c|}{6.0\%} & \multicolumn{2}{c|}{8.0\%} & \multicolumn{2}{c}{10.0\%} \\
        \midrule
        Clustering? & \ding{51} & \ding{55} & \ding{51} & \ding{55} & \ding{51} & \ding{55} \\
        \midrule
        Achievement Rate & 105.40\% & 43.81\% & 100.56\% & 39.11\% & 96.99\% & 40.72\% \\
        \midrule
        Participant Benefit & 0.2735 & 0.0958 & 0.4617 & 0.1847 & 0.6910 & 0.3363 \\
        \midrule
        \bottomrule
    \end{tabular}
    }
    \label{tab:clustering_micro}
\end{table}

\subsubsection{Necessity of Consumption Pattern Analysis}
We examine the impact of \emph{consumption pattern analysis} on assessing entity-level DR potential for flexible load regulation.
In this ablation, we fix the trade-off factor to \(\alpha = 10^3\) and compare DRP-FLR with its variant without pattern clustering in both the regional grid and the campus microgrid.
The results are summarized in TABLES~\ref{tab:clustering_regional} and ~\ref{tab:clustering_micro}.

In the regional-grid, clustering different load patterns for each entity yields an Achievement Rate deviation of \([-4.78\%,\,0.00\%]\) relative to the target interval \([1,\,1+\sigma]\),
which is substantially better than the \([5.17\%,\,90.87\%]\) deviation without clustering.
We attribute this over-regulation to the near-uniform distribution of load patterns among the enterprises.
When DRP-FLR estimates DR potential using a global historical mean, it tends to underestimate the actual adjustable range in many periods and compensates by selecting more entities, leading to excessive Achievement Rate.
Moreover, the over-regulation does not translate into a consistent improvement in Participant Benefit.
In particular, under \(\mathrm{Regulation\ Ratio}=10\%\), an over-regulation level of \(125.17\%\) leads to a \(24.76\%\) decrease in Participant Benefit.

In the campus microgrid, enabling pattern clustering yields an Achievement Rate deviation of \([-3.01\%,\,0.00\%]\), which is substantially better than the \([-60.89\%,\,-56.19\%]\) deviation without clustering.
This under-regulation is mainly caused by the fine-grained, submeter-level consumption patterns (e.g., floor-level HVAC), which exhibit a long-tail distribution.
When the regulation potential is estimated from global averages, DRP-FLR tends to overestimate in most periods, because rare but extreme high-reduction patterns disproportionately affect the average.
Consequently, an insufficient number of entities are selected, making the real reduction fail to meet the target demand and decreasing Participant Benefit.

\begin{figure*}[t]
  \centering
  \includegraphics[width=2.0\columnwidth]{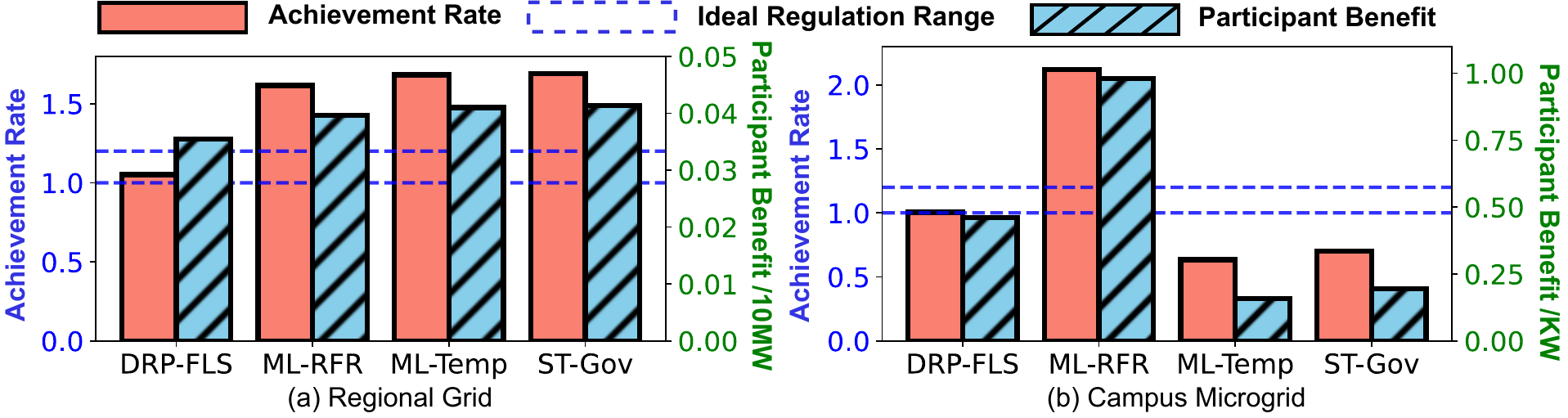}
  \caption{Comparison with Existing Demand Response Potential Estimation Methods.}
  \label{fig:contrast_DRP}
\end{figure*}

\begin{figure*}[t]
  \centering
  \includegraphics[width=2.0\columnwidth]{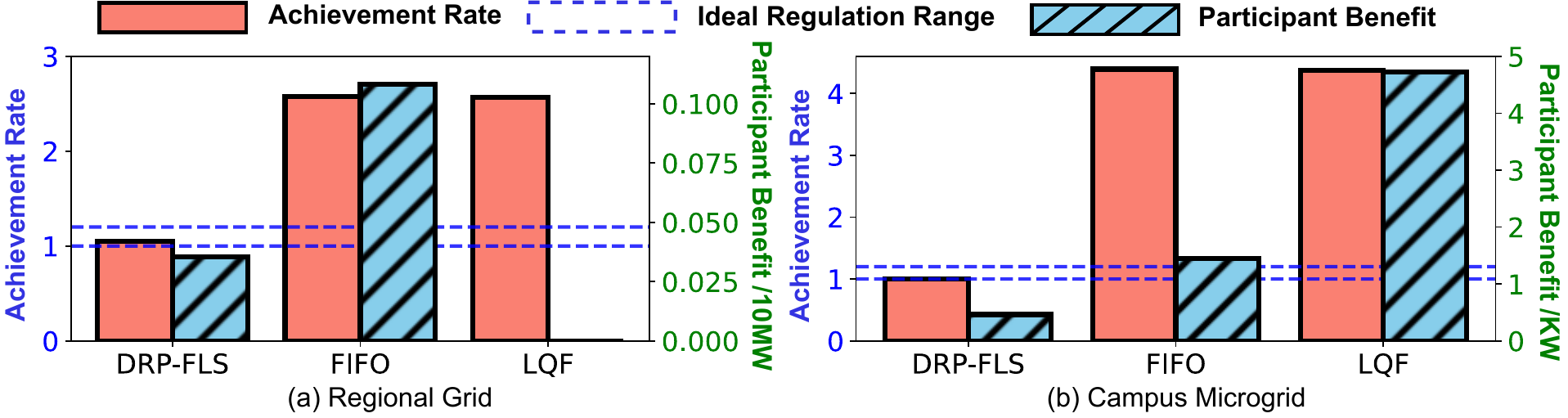}
  \caption{Comparison with Existing Naive Load Regulation Strategies.}
  \label{fig:contrast_Naive}
\end{figure*}

\subsection{Comparison Study}

\subsubsection{Different Potential Estimation Methods}

We compare DRP-FLR with representative load regulation strategies adopted by grid operators under DR mechanisms, where the key difference lies in how the regulation capability is estimated.
Li et al. propose an approach that requires both pre- and post-regulation load measurements for training~\cite{li2024online}, which are impractical to collect.
Therefore, we include three other practically applicable baselines in our comparison experiments:
\textbf{ML-RFR}~\cite{ruiz2024applications},
\textbf{ML-Temp}~\cite{muqtadir2025day}, and
\textbf{ST-Gov}~\cite{ruiz2019integration}.
We evaluate all methods under \(\mathrm{Regulation\ Ratio}=8\%\) and \(\alpha = 10^3\).
The results are summarized in Fig.~\ref{fig:contrast_DRP}.

In the regional grid, DRP-FLR achieves an Achievement Rate of \(105.15\%\), which falls within the ideal interval \([100\%,\,120\%]\) under \(\sigma=0.2\).
In contrast, the baseline methods exhibit pronounced over-regulation, exceeding the target by \(41.55\%\) to \(49.04\%\).
A key reason is that \textbf{ML-Temp} and \textbf{ST-Gov} do not distinguish multiple load patterns of each enterprise.
Due to the relatively uniform distribution of load patterns, their potential estimation tends to underestimate the real adjustable ranges.
Although \textbf{ML-RFR} partially captures variability via quantile regression forests, the overall load curves of each enterprise are relatively smooth (e.g., Fig.~\ref{fig:ex_pattern_regional_grid}).
This reduces quantile separability, also leading to underestimated regulation capacity and unnecessary participant selection.
While the baselines obtain higher \emph{total} Participant Benefit due to over-regulation, DRP-FLR improves the \emph{average} benefit by \(39.71\%\), better motivating participating entities.

In the campus microgrid, DRP-FLR achieves an Achievement Rate of \(100.56\%\), also within the desired interval.
In contrast, the baseline methods show substantial deviations ranging from \(-36.63\%\) to \(91.87\%\).
In the campus microgrid, submeter-level loads exhibit a pronounced long-tail distribution across appliances and subsystems.
Since \textbf{ML-Temp} and \textbf{ST-Gov} do not perform pattern classification, using global averages tends to overestimate the regulation capability, which reduces the number of selected entities and results in under-regulation.
Notably, \textbf{ML-RFR} substantially underestimates response potential, which in turn causes severe over-regulation, because extreme yet infrequent high consumption has little influence on the upper quantiles.
Besides, DRP-FLR increases the \emph{average} benefit by \(49.60\%\) compared with the baselines.

\subsubsection{Existing Naive Load Regulation Strategies}

In practice, many grid operators still rely on naive load regulation strategies without explicitly estimating participants’ DR potential.
Two widely adopted heuristics are \textbf{First-In-First-Out (FIFO)} and \textbf{Largest-Quantity-First (LQF)}.
FIFO selects participants according to their registration order until the demand is satisfied, whereas LQF prioritizes entities with larger contractual response quantities.
We compare DRP-FLR with these two naive baselines under the same setting as above, i.e., \(\mathrm{Regulation\ Ratio}=8\%\) and \(\alpha=10^3\), on both the regional grid and the microgrid.
The results are shown in Fig.~\ref{fig:contrast_Naive}.

Since FIFO and LQF lack any explicit assessment of entity-level response potential, both methods lead to severe over-regulation
In the regional grid, the deviation of Achievement Rate reaches \([136.93\%,\,137.51\%]\), while in the campus microgrid it further increases to \([217.59\%,\,219.35\%]\).
Such excessive regulation substantially violates the intended tolerance band \([1,\,1+\sigma]\) and can seriously compromise both grid supply-demand balance and economic feasibility (e.g., unnecessary dispatch costs and avoidable disruption to participants).

Notably, LQF even yields extremely low participant benefit in some cases.
This occurs because selecting entities with contractual response quantities does not guarantee that they are \emph{reward-effective} after settlement: large consumption does not necessarily imply large \emph{reducible} consumption.
The poor rewards for participants will undermine long-term incentives.

Overall, these results further highlight the importance of DRP-FLR: \emph{accurate DR potential estimation} is essential for achieving regulation targets without excessive over-regulation, while simultaneously sustaining participant incentives through reliable and fair benefit outcomes.


\section{Conclusion}
\label{concl}
The rapid advancement of AI and the large-scale integration of renewable energy exacerbate supply-demand imbalance in smart grids, motivating demand response mechanisms to enable peak shaving and valley filling.
However, existing operator-side load regulation strategies under demand response often ignore the diverse consumption patterns of various entities, resulting in biased potential estimation and severe over-regulation or under-regulation.
Therefore, we propose DRP-FLR. To our knowledge, it is the first load regulation method that achieves demand response potential assessment via consumption pattern analysis.
DRP-FLR jointly optimizes demand response potential, participant economic benefit, and renewable energy accommodation while ensuring supply-demand balance and economic feasibility of regulation plans.
Experiments on a real regional grid and a campus microgrid demonstrate that DRP-FLR can reduce regulation deviation by \(36.63\%\)-\(91.87\%\) and achieves an average \(44.66\%\) improvement in participant benefit.

\bibliographystyle{IEEEtran}
\bibliography{reference.bib}

\end{document}